# BAYESIAN NONPARAMETRIC ESTIMATORS DERIVED FROM CONDITIONAL GIBBS STRUCTURES

By Antonio Lijoi,[1] Igor Prünster[2] and Stephen G. Walker

*Università degli Studi di Pavia, Università degli Studi di Torino and University of Kent*

We consider discrete nonparametric priors which induce Gibbs-type exchangeable random partitions and investigate their posterior behavior in detail. In particular, we deduce conditional distributions and the corresponding Bayesian nonparametric estimators, which can be readily exploited for predicting various features of additional samples. The results provide useful tools for genomic applications where prediction of future outcomes is required.

**1. Introduction.** Random partitions and their associated probability distributions play an important role in a variety of research areas. In population genetics, for example, models for random partitions are useful in order to describe the allocation of a sample of $n$ genes into a number of distinct alleles. See, for example, [10, 33]. In machine learning theory, probabilistic models for linguistic applications (such as, e.g., speech and handwriting recognition, machine translation) are often based on a suitable clustering structure for a set of words. See, for example, [34, 35]. In Bayesian nonparametric inference, a discrete nonparametric prior is commonly employed in complex hierarchical mixture models and it induces an exchangeable random partition for the latent variables: this provides an effective tool for inferring on the clustering structure of the observations. Such an approach is due to [21] and has been extended in various directions. See, for example, [12, 13, 19, 22]. Other important areas of applications include storage problems, excursion theory,

Received May 2007; revised September 2007.
[1]Also affiliated with CNR-IMATI, Milan, Italy. Supported by MiUR, Grant 2006/134525.
[2]Also affiliated with Collegio Carlo Alberto and ICER, Turin, Italy. Supported by MiUR, Grant 2006/133449.

combinatorics and statistical physics. See the comprehensive and stimulating monograph by Pitman [29] and references therein.

An early and well-known model which describes the grouping of $n$ objects into $k$ distinct classes is due to [7] and leads to the Ewens sampling formula. The basic assumption is that individuals are sequentially sampled from an infinite set of different species and the proportion $\tilde{p}_i$ with which the $i$th species is present in the population is random. Then, if $(W_k)_{k\geq 1}$ is a sequence of independent and identically distributed random variables with $\text{Beta}(1,\theta)$ distribution, the random proportions are defined as

$$\tilde{p}_1 = W_1, \qquad \tilde{p}_j = W_j \prod_{k=1}^{j-1}(1-W_k) \qquad \forall j \geq 2. \tag{1}$$

Now, if $X_1, \ldots, X_n$ is a sample of $n$ individuals drawn from the population, set $\mathbf{M}_n := (M_{1,n}, \ldots, M_{n,n})$ where $M_{j,n}$ is the number of species represented $j$ times in the sample of size $n$. Hence, the distribution of $\mathbf{M}_n$ is supported by all those vectors $\mathbf{m}_n = (m_{1,n}, \ldots, m_{n,n})$ for which $\sum_{i=1}^{n} i m_{i,n} = n$. The Ewens sampling formula provides the probability distribution of the random vector $\mathbf{M}_n$ under (1) and it coincides with

$$\Pr[\mathbf{M}_n = \mathbf{m}_n] = \frac{n!}{(\theta)_n} \prod_{j=1}^{n} \frac{\theta^{m_{j,n}}}{j^{m_{j,n}} m_{j,n}!} \tag{2}$$

where $(\theta)_n = \theta(\theta+1)\cdots(\theta+n-1)$ for any $\theta > 0$. We also agree on setting $(\theta)_0 := 1$. See also [2] for a derivation of (2). Obviously, to the distribution of $\mathbf{M}_n$ there corresponds a distribution of the vector $(K_n, \mathbf{N}_n)$ where $K_n$ is the number of distinct species detected among the $n$ observations in the sample and $\mathbf{N}_n = (N_{1,n}, \ldots, N_{K_n,n})$ is the vector of frequencies with which each distinct species is observed. Such a correspondence is one-to-one and, conditional on $K_n$, the distribution of $\mathbf{N}_n$, is supported on the set $\Delta_{n,K_n} := \{(n_1, \ldots, n_{K_n}) \in \{1,\ldots,n\}^{K_n} : \sum_{j=1}^{K_n} n_j = n\}$. In particular, for the Ewens sampling formula (2) there corresponds the probability distribution

$$\Pr[K_n = k, \mathbf{N}_n = (n_1, \ldots, n_{K_n})] = \frac{\theta^k}{(\theta)_n} \prod_{j=1}^{k} (n_j - 1)! \tag{3}$$

for any $k \in \{1,\ldots,n\}$ and $(n_1,\ldots,n_k) \in \Delta_{n,k}$. The parameter $\theta$, in genetic applications, is interpreted as the mutation rate of each gene into new allelic types. Formula (3) has a further interesting combinatorial interpretation. If $\theta$ is a positive integer, then $\theta^k \prod_{j=1}^{k}(n_j - 1)!$ is the number of colored permutations of $\{1,\ldots,n\}$ into $k$ cycles with respective lengths $n_1,\ldots,n_k$, each cycle being labeled by any of the $\theta$ available colors. Accordingly, (3) is the probability distribution of a random permutation with colored cycles. See [3, 29] for exhaustive accounts on the Ewens sampling formula.

The distribution of the vector $(K_n, \mathbf{N}_n)$ takes on the name of *exchangeable partition probability function* (EPPF), a notion introduced by Pitman in [26] and further studied in a series of subsequent papers; see [29] and references therein. The main object of investigation of the present paper is a family of EPPFs, introduced and thoroughly investigated in [9], which generalize the Ewens sampling scheme. Our aim is to establish distributional properties of such EPPFs which allow, given a sample, to make predictions according to a Bayesian nonparametric procedure. The concrete motivation for this study is provided by the straightforward applicability of the results to inference in genetic experiments. As a matter of fact, an important setting where our findings can be usefully applied relates to gene detection in expressed sequence tags (EST) experiments. ESTs are produced by sequencing randomly selected cDNA clones from a cDNA library. Given an initial EST data set of size $n$, one is interested in the prediction of the outcomes of further sampling from the library. For instance, interest lies in the estimation of the number of new unique genes in a possible additional sample of size $m$: nonparametric frequentist estimators, however, yield completely unstable estimates when $m > 2n$. See [25] for a discussion of this phenomenon. In contrast, for the corresponding Bayesian nonparametric estimators proposed in [20], and based on Gibbs partitions, the relative dimension of $m$ with respect to $n$ is not an issue. Indeed, we will show that the EPPF, whenever analytically available, yields straightforward and coherent answers to this and other related prediction problems.

In Section 2 we recall the concepts of exchangeable random partition and EPPF and the definition of the class of exchangeable Gibbs random partitions. In Section 3 we derive distributional results for the corresponding EPPFs conditionally on a sample: we obtain expressions for the predictive distribution of future observations given the past, then focus on the probability distribution of the random partition restricted to those observations yielding new distinct species in the future sample and, finally, face the problem of determining the probability that specific observed species will not appear in the future sample. In Section 4 we illustrate how our results can be applied in the context of EST analysis of cDNA libraries. The Appendix contains a short review of generalized factorial coefficients and the proofs.

**2. Exchangeable Gibbs random partitions.** A random partition of the set of natural numbers $\mathbb{N}$ is defined as a *consistent* sequence $\Pi = \{\Pi_n\}_{n=1}^{\infty}$ of random elements, with $\Pi_n$ taking values in the set of all partitions of $[n] := \{1, \ldots, n\}$ into some number of disjoint nonempty blocks. Consistency in this setting implies that each $\Pi_n$ is obtained from $\Pi_{n+1}$ by discarding the integer $n+1$. A random partition $\Pi$ is *exchangeable* if, for each $n$, the probability distribution of $\Pi_n$ is invariant under all permutations of $(1, \ldots, n)$. To be more precise, let $\{A_j\}_{j=1}^k$ denote a partition of the set $[n]$, and let the $A_j$'s be indexed by $[k]$ in order of their least elements. In order to describe the property of exchangeability for $\Pi$ let us introduce a sequence of functions $\Pi_k^{(n)} : \Delta_{n,k} \to \mathrm{R}^+$ such that:

(i) $\Pi_1^{(1)}(1) = 1$;
(ii) for any $(n_1, \ldots, n_k) \in \Delta_{n,k}$, $k \in \{1, \ldots, n\}$ and $n \geq 1$ one has

$$\Pi_k^{(n)}(n_1, \ldots, n_k) = \Pi_k^{(n)}(n_{\rho(1)}, \ldots, n_{\rho(k)})$$

where $\rho$ is an arbitrary permutation of the indices $(1, \ldots, k)$;

(iii) for any $(n_1, \ldots, n_k) \in \Delta_{n,k}$, $k \in \{1, \ldots, n\}$ and $n \geq 1$ the following addition rule holds true:

$$\begin{aligned}(4) \quad & \Pi_k^{(n)}(n_1, \ldots, n_k) \\ & = \sum_{j=1}^k \Pi_k^{(n+1)}(n_1, \ldots, n_j+1, \ldots, n_k) + \Pi_{k+1}^{(n+1)}(n_1, \ldots, n_k, 1).\end{aligned}$$

A function $\Pi_k^{(n)}$ with these properties is known as an *exchangeable partition probability function* (EPPF) and it uniquely determines the probability law of an exchangeable random partition according to the equality

$$\mathbb{P}(\Pi_n = \{A_1, \ldots, A_k\}) = \Pi_k^{(n)}(|A_1|, \ldots, |A_k|),$$

where $|A|$ stands for the cardinality of set $A$. A first treatment of this concept can be found in [26], and a recent exhaustive account on exchangeable random partitions is provided in [29]. The above-mentioned Ewens sampling formula corresponds to the EPPF of the Dirichlet process [8] as described in (3) and it has found many interesting applications, for instance, in Bayesian nonparametrics and in population genetics. Another noteworthy example is represented by Pitman's sampling formula which corresponds to an EPPF of the form

$$(5) \quad \Pi_k^{(n)}(n_1, \ldots, n_k) = \frac{\prod_{i=1}^{k-1}(\theta + i\sigma)}{(\theta + 1)_{n-1}} \prod_{j=1}^k (1 - \sigma)_{n_j - 1},$$

where $\theta > -\sigma$ and $\sigma \in (0, 1)$ or $\sigma < 0$ and $\theta = \nu|\sigma|$ for some positive integer $\nu$. See [26]. This can also be seen as the probability distribution induced by the species sampling model $\tilde{P}(\cdot) = \sum_{j=1}^\infty \tilde{p}_j \delta_{X_j}$ where the $X_j$'s are independent and identically distributed from some nonatomic distribution $H$ and the weights $\tilde{p}_j$ are constructed via a stick-breaking procedure as in (1) the only difference being, now, that $W_j \sim \mathrm{Beta}(1 - \sigma, \theta + j\sigma)$ for any $j \geq 1$. We also agree that $W_j \sim \mathrm{Beta}(1 - \sigma, 0)$ implies that $W_j = 1$ almost surely. The random probability $\tilde{P}$ is termed the two parameter Poisson–Dirichlet process. See [27, 29].

Another interesting example of EPPF arises from the normalization of a generalized gamma process, as defined in [4], and leads to

$$\Pi_k^{(n)}(n_1,\ldots,n_k) \tag{6}$$
$$= \frac{\sigma^{k-1}e^\beta \prod_{j=1}^k (1-\sigma)_{n_j-1}}{\Gamma(n)} \sum_{i=0}^{n-1} \binom{n-1}{i}(-1)^i \beta^{i/\sigma}\Gamma\left(k-\frac{i}{\sigma};\beta\right)$$

where $\beta > 0$ and $\Gamma(a;x) := \int_x^\infty s^{a-1}e^{-s}\,ds$ is, for any $x>0$, the incomplete gamma function. See [14, 28] and [19] for an application of the corresponding random discrete distribution in the context of mixture modeling. For general results concerning random probability measures derived via normalization procedures see [15, 16, 17, 28, 31].

The examples we have briefly illustrated so far share a common structure. Indeed, one may note that each EPPF in (3), (5) and (6) arises as a product of two factors: the first one depends only on $(n,k)$ and the second one depends on the frequencies $(n_1,\ldots,n_k)$ via the product $\prod_{j=1}^k (1-\sigma)_{n_j-1}$. This structure is the main ingredient for defining a general family of exchangeable random partitions, namely the *Gibbs-type random partitions*.

DEFINITION 1 ([9]). An exchangeable random partition $\Pi$ of the set of natural numbers is said to be of *Gibbs form* if, for all $1 \leq k \leq n$ and for any $(n_1,\ldots,n_k)$ in $\Delta_{n,k}$, the EPPF of $\Pi$ can be represented as

$$\Pi_k^{(n)}(n_1,\ldots,n_k) = V_{n,k}\prod_{j=1}^k (1-\sigma)_{n_j-1}, \tag{7}$$

for some $\sigma \in [0,1)$.

It is worth noting that the previous definition holds also for negative values of $\sigma$. See [9]. According to Definition 1, an exchangeable Gibbs-type random partition is completely specified once the $V_{n,k}$'s have been assigned. As shown in [9], if a set of nonnegative weights $\mathscr{V} := \{V_{n,k} : k=1,\ldots,n; n \geq 1\}$ solves the forward recursive equations

$$V_{n,k} = (n-\sigma k)V_{n+1,k} + V_{n+1,k+1}, \tag{8}$$

then $\mathscr{V}$ identifies the EPPF of a Gibbs-type random partition. Hence, for infinite exchangeable sequences of random partitions, the above recursion might provide a constructive approach in order to determine Gibbs-type random partitions. In [9], Theorem 12, one can find a complete description of the extreme points of $\mathscr{V}$. With reference to the previously illustrated examples, the corresponding set of weights $\mathscr{V}$ are immediately identified from (3), (5) and (6), respectively. Recently, [11] have investigated the dependence of the distribution of the frequencies of the clusters of a Gibbs-type partition on their least elements and have extended some of the results contained in [10] relating to the Ewens sampling formula.

Finally, note that Definition 1 directly involves infinite sequences $\Pi = \{\Pi_n\}$ of exchangeable random partitions. One can, however, confine oneself to considering just a finite sequence of partitions $\Pi = \{\Pi_n\}_{n=1}^N$ for some integer $N \geq 1$. In this case, we say that $\Pi$ is a *finite Gibbs random partition* if it is characterized by an EPPF of the form (7), for any $k \in \{1, \ldots, n\}$ and $n \in \{1, \ldots, N\}$. Note that in this case, the addition rule (4) defining the EPPF holds true for $n \in \{1, \ldots, N-1\}$.

**3. Conditional structures of Gibbs-type random partitions.** The main goal we are pursuing in the present paper consists in investigating some conditional structures that emerge when the observations are sampled according to a Gibbs-type random partition with a view to deriving Bayesian nonparametric estimators for quantities of interest. The issue we address consists in evaluating, conditionally on the partition of a basic sample of size $n$, the probability of sampling, if $m$ draws, a certain number of observations yielding new partition groups with specified frequencies. Such a quantity can be useful in a variety of applications, some of which we highlight in Section 4. Resorting to the notation set forth in the Section 2, we study distributional properties of the partition of the set of integers $\{n+1, \ldots, n+m\}$, given $[n]$ has been partitioned into $j$ classes with respective frequencies $(n_1, \ldots, n_j)$. A few quantities, analogous to those describing the partition structure of $[n]$, need to be introduced in advance. We let $K_m^{(n)} = K_{m+n} - K_n$ stand for the number of new partition sets $C_1, \ldots, C_{K_m^{(n)}}$ generated by the additional $m$-sample $X_{n+1}, \ldots, X_{n+m}$. Furthermore, if $C := \bigcup_{i=1}^{K_m^{(n)}} C_i$ whenever $K_m^{(n)} \geq 1$ and $C \equiv \varnothing$ if $K_m^{(n)} = 0$, we set $L_m^{(n)} := \mathrm{card}(\{X_{n+1}, \ldots, X_{n+m}\} \cap C)$ as the number of observations belonging to the new clusters $C_i$. It is obvious that $L_m^{(n)} \in \{0, 1, \ldots, m\}$ and that $m - L_m^{(n)}$ observations belong to the sets defining the partition of the original $n$ observations. According to this, if $\mathbf{S}_{L_m^{(n)}} = (S_{1, L_m^{(n)}}, \ldots, S_{K_m^{(n)}, L_m^{(n)}})$ then the distribution of $\mathbf{S}_{L_m^{(n)}}$, conditional on $L_m^{(n)} = s$, is supported by all vectors $(s_1, \ldots, s_{K_m^{(n)}})$ of positive integers such that $\sum_{i=1}^{K_m^{(n)}} s_i = s$. The remaining $m - L_m^{(n)}$ observations are allocated to the "old" $K_n$ groups with vector of nonnegative frequencies $\mathbf{R}_n = (R_1, \ldots, R_{K_n})$ such that $\sum_{i=1}^{K_n} R_i = m - L_m^{(n)}$. Throughout we also assume that all random quantities are defined on a common probability space $(\Omega, \mathscr{F}, \mathbb{P})$.

PROPOSITION 1. *Suppose that $\Pi = \{\Pi_n\}_{n=1}^\infty$ is a Gibbs-type exchangeable random partition with weights $V_{n,k}$ and parameter $\sigma \in [0, 1)$. Then, the* *joint distribution of $K_m^{(n)}$, $L_m^{(n)}$ and $\mathbf{S}_{L_m^{(n)}}$, given $K_n$ and $\mathbf{N}_n$, is of the form*

$$\mathbb{P}(K_m^{(n)} = k, L_m^{(n)} = s, \mathbf{S}_{L_m^{(n)}} = (s_1, \ldots, s_{K_m^{(n)}}) | K_n = j, \mathbf{N}_n = (n_1, \ldots, n_{K_n}))$$

(9) $\quad = \mathbb{P}(K_m^{(n)} = k, L_m^{(n)} = s, \mathbf{S}_{L_m^{(n)}} = (s_1, \ldots, s_{K_m^{(n)}}) | K_n = j)$

$$= \frac{V_{n+m,j+k}}{V_{n,j}} \binom{m}{s} (n - j\sigma)_{m-s} \prod_{i=1}^{k} (1 - \sigma)_{s_i - 1}.$$

*Hence, the number $K_n$ of partition sets in the basic $n$ sample is sufficient for predicting:* (i) *the number of sets into which $\{n+1, \ldots, n+m\}$ is partitioned,* (ii) *the number of points from the subsequent $m$ sample that belong to the new sets of the partition of $[n+m]$ and* (iii) *the frequencies in each of these new groups.*

By marginalizing the conditional distribution in (9) with respect to $\mathbf{S}_{L_m^{(n)}}$ and, then, with respect to $K_m^{(n)}$ one obtains the conditional distribution for the number of new groups and the number of observations belonging to these new groups and the distribution of $L_m^{(n)}$, respectively. These marginalizations yield results in terms of *generalized Stirling numbers* or *generalized factorial coefficients*, denoted as $\mathscr{C}(s, k, \sigma)$ and whose representation is given in (37).

COROLLARY 1. *The joint distribution of $K_m^{(n)}$ and $L_m^{(n)}$, given $K_n$, can be expressed as*

(10)
$$\mathbb{P}(K_m^{(n)} = k, L_m^{(n)} = s | K_n = j)$$
$$= \frac{V_{n+m,j+k}}{V_{n,j}} \binom{m}{s} (n - j\sigma)_{m-s} \frac{\mathscr{C}(s, k, \sigma)}{\sigma^k}$$

*for $k \leq s = 0, \ldots, m$ and the conditional distribution of $L_m^{(n)}$ is of the form*

(11) $\quad \mathbb{P}(L_m^{(n)} = s | K_n = j) = \binom{m}{s} (n - j\sigma)_{m-s} \sum_{k=0}^{s} \frac{V_{n+m,j+k}}{V_{n,j}} \frac{\mathscr{C}(s, k, \sigma)}{\sigma^k}$

*for $s = 0, \ldots, m$.*

From (10) and (11) one can also deduce other explicit forms for conditional distributions of interest. For example, the distribution of the number of observations in the new $m$-sample which lie in new partition sets, given the number of groups present in the basic $n$-sample and the number of new clusters $K_m^{(n)}$, is of the form

(12) $\quad \mathbb{P}(L_m^{(n)} = s | K_m^{(n)} = k, K_n = j) = \dfrac{\binom{m}{s}(n - j\sigma)_{m-s} \mathscr{C}(s, k, \sigma)}{\mathscr{C}(m, k; \sigma, -n + j\sigma)}$

for $s = k, \ldots, m$, where $\mathscr{C}(n,k;\sigma,\gamma)$ is a noncentral generalized factorial coefficient representable as in (39). It is worth noting that the previous expression does not depend on the particular Gibbs prior it is derived from: interestingly, Gibbs-type random partitions share the same conditional structures once $K_m^{(n)}$ and $K_n$ are fixed. This finding is reminiscent of a result in [9] where the authors show that $K_n$ is sufficient for the Gibbs random partition of the first $n$ integers meaning that the conditional distribution of the partition of $[n]$ given $K_n$ does not depend on the weights $V_{n,k}$. On the other hand, the conditional distribution of $K_m^{(n)}$, given $L_m^{(n)}$ and $K_n$, is of the form

$$(13) \qquad \mathbb{P}(K_m^{(n)} = k | L_m^{(n)} = s, K_n = j) = \frac{V_{n+m,j+k}\mathscr{C}(s,k,\sigma)/\sigma^k}{\sum_{l=0}^{s} V_{n+m,j+l}\mathscr{C}(s,l,\sigma)/\sigma^l}$$

for any $k \in \{0, \ldots, s\}$. Moreover, the Bayes estimator (under quadratic loss function) for the expected number of new clusters, proposed in [20], is easily recovered from (10) as

$$(14) \qquad \mathbb{E}(K_m^{(n)} | K_n = j) = \sum_{k=0}^{m} k \frac{V_{n+m,j+k}}{V_{n,j}} \frac{\mathscr{C}(m,k;\sigma,-n+j\sigma)}{\sigma^k}.$$

Often interest relies also in determining an estimator for the number of observations in the subsequent $m$-sample that will belong to new species. For instance, in genomic applications this can be seen as a better measure of redundancy of a certain library. For this purpose, one can resort to (11) and the corresponding Bayes estimator is given by

$$(15) \quad \mathbb{E}(L_m^{(n)} | K_n = j) = \sum_{s=0}^{m} s \binom{m}{s} (n - j\sigma)_{m-s} \sum_{k=0}^{s} \frac{V_{n+m,j+k}}{V_{n,j}} \frac{\mathscr{C}(s,k,\sigma)}{\sigma^k}.$$

Then, $\mathbb{E}(L_m^{(n)} | K_n = j)/m$ is the expected proportion of genes in the new sample which do not coincide with previously observed ones. The expression in (15) admits a noteworthy simplification as outlined in the following proposition: indeed, the Bayes estimator is $m$ times the probability that the $(n+1)$th draw yields a new cluster, given that $j$ distinct clusters are generated by the first $n$ observations.

PROPOSITION 2. *For any $j \in \{1, \ldots, n\}$ and $m \geq 1$ one has*

$$(16) \qquad \mathbb{E}(L_m^{(n)} | K_n = j) = m \frac{V_{n+1,j+1}}{V_{n,j}}.$$

All the previous expressions are easily available for the three examples we have mentioned in Section 2. We first focus our attention on the Dirichlet process which represents the most well-known case. Indeed, from (3) one finds out that $V_{n,k} = \theta^k/(\theta)_n$ and, for instance, (9) reduces to

(17)
$$\mathbb{P}(L_m^{(n)} = s, K_m^{(n)} = k, \mathbf{S}_{L_m^{(n)}} = (s_1, \ldots, s_{K_m^{(n)}}) | K_n = j)$$
$$= \frac{\theta^k}{(\theta + n)_m} \binom{m}{s} (n)_{m-s} \prod_{i=1}^{k} (s_i - 1)!.$$

Note that simple algebra leads to rewrite the above expression as

$$\binom{m}{s} \left(1 - \frac{n}{\theta + n}\right)^k \left\{ \frac{(\theta + n)^k}{(\theta + n)_s} \prod_{i=1}^{k} (s_i - 1)! \right\} \frac{(n)_{m-s}}{(\theta + n + s)_{m-s}}$$
$$= \binom{m}{s} p_\theta(n, m, k, s, \mathbf{s}_k)$$

where it can be immediately seen that the term in curly brackets on the left-hand side is the sampling formula in (9) with $\theta + n$ in the place of $\theta$ being the total mass parameter of the Dirichlet process conditioned on a sample of size $n$. Hence, the quantity $p_\theta(n, m, k, s, \mathbf{s}_k)$ can be interpreted as the probability of drawing, conditional on the $n$ past observations, a specific sample of size $m$ of which $s$ belong to the new $k$ groups of the partition with vector of frequencies $\mathbf{s}_k = (s_1, \ldots, s_k)$ and the other $m - s$ coincide with any of the conditioning $n$ observations. On the other hand, recall that $\lim_{\sigma \to 0} \frac{\mathscr{C}(n,k,\sigma)}{\sigma^k} = |\mathfrak{s}(n,k)|$ where $\mathfrak{s}(\cdot, \cdot)$ stands for the Stirling number of the first kind. This allows to determine the expressions appearing in (10) and (11). Indeed, one has

$$\mathbb{P}(K_m^{(n)} = k, L_m^{(n)} = s | K_n = j) = \binom{m}{s} \frac{\theta^k (n)_{m-s}}{(\theta + n)_m} |\mathfrak{s}(n, k)|$$

and, using the definition of the signless Stirling number of the first kind according to which

(18)
$$\sum_{i=0}^{s} \theta^i |\mathfrak{s}(s, i)| = (\theta)_s$$

(see, e.g., [6], page 2536), one has

$$\mathbb{P}(L_m^{(n)} = s | K_n = j) = \binom{m}{s} \frac{(n)_{m-s}}{(\theta + n + s)_{m-s}} \frac{(\theta)_s}{(\theta + n)_s} = \binom{m}{s} q_\theta(n, m, s)$$

where it is apparent that $q_\theta(n, m, s)$ is the probability, conditional on a sample of size $n$, of observing a specific $m$-sample containing $s$ elements not contained in the conditioning $n$-sample.

REMARK 1. It is important to note that the conditional structure of the Dirichlet process does not depend on $K_n$: it only depends on the size of the basic sample $n$. This is, indeed, a characterizing property of the Dirichlet process as shown in [36]. Such a property simplifies the mathematical expressions but represents a serious drawback for applications. Indeed, it is reasonable to expect that $K_n$ influences prediction of the clustering structure of future observations: the larger $K_n$ the more new clusters $K_m^{(n)}$ and the more observations belonging to these new clusters $L_m^{(n)}$ one would expect. This is the reason which explains the interest in a more general family of partition distributions such as those of Gibbs-type for which prediction depends on $K_n$. Finally, it is worth recalling that the Dirichlet process can be seen as a two parameter Poisson–Dirichlet process with parameter $(\theta, 0)$. Hence, when we deal with the Poisson–Dirichlet process in the sequel, the Dirichlet process case can be recovered by letting $\sigma \to 0$.

REMARK 2. All the quantities described up to now, and developed in the next subsections, depend on the analysis of the conditional structure of a Gibbs-type random partition. Investigation of the conditional structure for the sequence of blocks $(K_n)_{n \geq 1}$ is pursued in [9] where the authors do consider the conditional distribution of the number of groups in the partition of $[n]$, given the number of blocks in which $[n+m]$ is partitioned. In our setting, where prediction is the main focus, we are more interested in evaluating conditional probabilities (or expectations) for the partition of future observations given the partition structure of past observations. And we also consider other relevant quantities, besides the number of groups. It might be that starting from the conditional characterizations provided by [9] one can derive formulae analogous to those we are now going to establish, but we find our approach more direct and particularly suited to the specific prediction problems we have in mind.

3.1. *The process generating new clusters.* We are now going to consider an important quantity which describes the partition structure of observations generating new groups in a further sampling procedure, conditional on the partition generated by the first $n$ observations. In particular we are able to point out a sort of reproducibility of the Gibbs structure as established by the following proposition.

PROPOSITION 3. *Let $\Pi = \{\Pi_n\}_{n=1}^{\infty}$ be a Gibbs-type random exchangeable partition whose EPPF is characterized by the set of weights $\{V_{n,k} : k = 1, \ldots, n; n \geq 1\}$ and by the parameter $\sigma \in (0,1)$. Then*

$$\mathbb{P}(K_m^{(n)} = k, \mathbf{S}_{L_m^{(n)}} = (s_1, \ldots, s_{K_m^{(n)}}) | L_m^{(n)} = s, K_n = j, \mathbf{N}_n = (n_1, \ldots, n_j)) \quad (19)$$

$$= \frac{V_{n+m,j+k}}{\sum_{i=0}^{s} V_{n+m,j+i}\mathscr{C}(s,i,\sigma)/\sigma^i} \prod_{i=1}^{k}(1-\sigma)_{s_i-1}$$

*for any $s \in \{1,\ldots,m\}$, $k \in \{1,\ldots,s\}$, $j \in \{1,\ldots,n\}$, $(n_1,\ldots,n_j) \in \Delta_{n,j}$ and $(s_1,\ldots,s_k) \in \Delta_{s,k}$. Consequently the partition of the observations which belong to the new partition sets is, conditional on the basic sample of size $n$, a finite Gibbs-type random partition with weights $\{V_{s,k}(m,n,j) : s = 1,\ldots,m; k = 1,\ldots,s\}$ defined by*

$$(20) \quad V_{s,k}(m,n,j) = \frac{V_{n+m,j+k}}{\sum_{i=0}^{s} V_{n+m,j+i}\frac{\mathscr{C}(s,i,\sigma)}{\sigma^i}}$$

*and with parameter $\sigma \in [0,1)$.*

Note from (19), again, that

$$\mathbb{P}(K_m^{(n)} = k, \mathbf{S}_{L_m^{(n)}} = (s_1,\ldots,s_{K_m^{(n)}})|L_m^{(n)} = s, K_n = j, \mathbf{N}_n = (n_1,\ldots,n_j))$$
$$= \mathbb{P}(K_m^{(n)} = k, \mathbf{S}_{L_m^{(n)}} = (s_1,\ldots,s_{K_m^{(n)}})|L_m^{(n)} = s, K_n = j).$$

The finiteness of the random partition described by (19) is obvious, since it takes values on the space of all partitions of $[s]$, with $1 \leq s \leq m$. Moreover, the particular structure featured by the conditional distribution in (19) motivates the following definition.

DEFINITION 2. *The conditional probability distribution*

$$(21) \quad \begin{aligned} &\tilde{\Pi}_k^{(s)}(s_1,\ldots,s_k; m,n,j) \\ &\quad := \mathbb{P}(K_m^{(n)} = k, \mathbf{S}_{L_m^{(n)}} = (s_1,\ldots,s_{K_m^{(n)}})|L_m^{(n)} = s, K_n = j), \end{aligned}$$

*with $1 \leq s \leq m$ and $1 \leq k \leq s$, is termed conditional EPPF.*

Hence, the probability distribution in (19) is a conditional EPPF giving rise to a finite Gibbs-type random partition. Even if the structure of $\tilde{\Pi}_k^{(s)}(s_1,\ldots,s_k; m,n,j)$ is quite general, one might wonder whether it is possible to provide more information about its $V_{s,k}(m,n,j)$ weights in some particular cases. For example, it would be interesting to ascertain when $V_{s,k}(m,n,j)$ does not depend on $m$ and $n$, so that $\tilde{\Pi}_k^{(s)}(s_1,\ldots,s_k; m,n,j) = \tilde{\Pi}_k^{(s)}(s_1,\ldots,s_k; j)$, which means that the conditional EPPF is that corresponding to an infinite Gibbs partition. This leads us to state the following:

COROLLARY 2. *The conditional EPPF $\tilde{\Pi}_k^{(s)}(s_1,\ldots,s_k; m,n,j)$ does not depend directly on $m$ and $n$ if and only if it is determined from a two-parameter Poisson–Dirichlet random partition.*

Having the conditional EPPF $\tilde{\Pi}_k^{(s)}$ at hand, one can compute some other interesting conditional distributions in a straightforward way. For example, if one combines the expression for $\tilde{\Pi}_k^{(s)}$ with Corollary 1 it is immediate to check that

$$\mathbb{P}(\mathbf{S}_{L_m^{(n)}} = (s_1, \ldots, s_{K_{m^{(n)}}}) | K_m^{(n)} = k, L_m^{(n)} = s, K_n = j)$$
$$= \frac{\sigma^k}{\mathscr{C}(s,k,\sigma)} \prod_{i=1}^{k} (1-\sigma)_{s_i - 1}$$

is an expression for the conditional distribution of detecting a particular configuration $(s_1, \ldots, s_k)$ for the observations belonging to the new partition sets, given the number of new sets, the number of observations falling into these sets and the basic $n$-sample.

All the sampling *formulae* we have deduced so far have important applications in Bayesian nonparametrics and population genetics. In Bayesian nonparametrics, random discrete distributions are commonly employed in order to define a clustering structure either at the level of the observations or at the level of the latent variables in a complex hierarchical model. In particular any EPPF corresponds to some random discrete distribution and it represents, together with all the expressions for the conditional distributions we have obtained, a useful tool for specifying prior opinions on the clustering of the data. In population genetics, the concept of conditional EPPF can be seen as follows. Given a sample of size $n$ containing $j$ distinct species with absolute frequencies $n_1, \ldots, n_j$, a new sample of size $m$ is to be drawn. Given that $s$ of the $m$ observations contribute to generating newly observed species, that is, they belong to new distinct clusters, one might be interested in evaluating the probability that the $s$ observations are grouped into $k$ clusters with respective frequencies $s_1, \ldots, s_k$. The answer to such a question is provided by a conditional EPPF. The other distributions, discussed previously, provide a wide range of sampling *formulae* which answer similar types of problems. In the following subsection we focus attention on some noteworthy particular cases, namely the Poisson–Dirichlet distribution, the two-parameter Poisson–Dirichlet distribution and the generalized gamma partition distribution.

3.2. *Illustrative examples.* We start our illustrations by considering the two-parameter Poisson–Dirichlet process due to [26]. The EPPF of this process is also known as *Pitman sampling formula*. Basing upon Proposition 1, one has

(22)
$$\mathbb{P}(K_m^{(n)} = k, L_m^{(n)} = s, \mathbf{S}_{L_m^{(n)}} = (s_1, \ldots, s_{K_m^{(n)}}) | K_n = j)$$
$$= \frac{\prod_{i=0}^{k-1}(\theta + j\sigma + i\sigma)}{(\theta + n)_m} \binom{m}{s} (n - j\sigma)_{m-s} \prod_{i=1}^{k} (1-\sigma)_{s_i - 1}$$

and it is possible to derive explicit expressions for all the sampling *formulae* set forth in Section 2. First note that from properties of generalized factorial coefficients, one has

$$\sum_{k=0}^{s} V_{n+m,j+k} \frac{\mathscr{C}(s,k,\sigma)}{\sigma^k} = \frac{\sigma^j}{(\theta+n)_m} \sum_{k=0}^{\sigma} \left(\frac{\theta}{\sigma}\right)_{j+k} \mathscr{C}(s,k,\sigma)$$

$$= \frac{\prod_{i=0}^{j-1}(\theta+i\sigma)}{(\theta+n)_m} \sum_{k=0}^{\sigma} \left(\frac{\theta}{\sigma}+j\right)_k \mathscr{C}(s,k,\sigma)$$

$$= \frac{\prod_{i=0}^{j-1}(\theta+i\sigma)}{(\theta+n)_m} (\theta+j\sigma)_s$$

$$= V_{n+m,j}(\theta+j\sigma)_s.$$

According to this equality, from (11) one has

$$(23) \qquad \mathbb{P}[L_m^{(n)} = s | K_n = j] = \frac{1}{(\theta+n)_m} \binom{m}{s} (n-j\sigma)_{m-s} (\theta+j\sigma)_s.$$

Now, (23) yields an estimate for the expected number of observations which do not coincide with the previously observed ones which, by virtue of Proposition 2, coincides with

$$(24) \qquad \mathbb{E}[L_m^{(n)} | K_n = j] = \frac{m(\theta+j\sigma)}{\theta+n}.$$

Consider now the conditional EPPF in (19), which is associated to the process generating the new clusters as explained in Section 3.1. We know by Corollary 2 that in the two-parameter Poisson–Dirichlet case the $V_{s,k}(m,n,j)$ weights do not depend on $m$ and $n$. Their specific form is easily seen to be

$$V_{s,k}(m,n,j) = \frac{V_{n+m,j+k}}{\sum_{i=0}^{s} V_{n+m,j+i}\sigma^{-i}\mathscr{C}(s,i,\sigma)} = \frac{V_{n+m,j+k}}{V_{n+m,j}(\theta+j\sigma)_s}$$

$$= \frac{\prod_{i=j}^{j+k-1}(\theta+i\sigma)}{(\theta+j\sigma)_s} = \frac{\prod_{i=0}^{k-1}(\theta+j\sigma+i\sigma)}{(\theta+j\sigma)_s}$$

with the proviso that $\prod_{i=s}^{s-1}(\theta+j\sigma+i) \equiv 1$. Hence, the conditional Pitman sampling formula is given by

$$(25) \qquad \tilde{\Pi}_k^{(s)}(s_1,\ldots,s_k;m,n,j) = \frac{\prod_{l=0}^{k-1}(\theta+j\sigma+l\sigma)}{(\theta+j\sigma)_s} \prod_{i=1}^{k}(1-\sigma)_{s_i-1}.$$

Now set $\theta' = \theta + j\sigma$ and note that the conditional EPPF of a Poisson–Dirichlet process with parameter $(\theta,\sigma)$ is again a Poisson–Dirichlet process with an updated parameter $(\theta',\sigma)$. This can be seen as a quasi-conjugacy of the two-parameter Poisson–Dirichlet process, where by quasi-conjugacy we mean that the process generating the new observations is of the same form as the prior process with updated parameters. Hence, at this stage one can equivalently re-express Corollary 2 above in a language quite familiar in Bayesian nonparametrics as follows:

COROLLARY 2'. *The only quasi-conjugate Gibbs-type prior is the two-parameter Poisson–Dirichlet process.*

Note that the quasi-conjugacy of the two-parameter Poisson–Dirichlet process was first shown in Pitman ([28], Corollary 20) by means of different techniques, whereas the characterization as the only quasi-conjugate Gibbs prior is new. With the Poisson–Dirichlet process with parameter $(\theta, 0)$, that is, the Dirichlet process prior with parameter measure having total mass $\theta > 0$, some useful simplifications occur. For example, the conditional EPPF is

$$(26) \quad \tilde{\Pi}_k^{(s)}(s_1, \ldots, s_k; m, n, j) = \frac{\theta^k \prod_{i=1}^k (s_i - 1)!}{\sum_{i=0}^s \theta^i |\mathfrak{s}(s, i)|} = \frac{\theta^k}{(\theta)_s} \prod_{i=1}^k (s_i - 1)!$$

which replicates the unconditional form of the EPPF in (3) and, as expected, does not depend on $(m, n, j)$. From a Bayesian nonparametric perspective, this is not surprising given the *conjugacy* of the Dirichlet process (see [8]). Indeed, this is just a reformulation, in a different context, of the fact that given a sample from the Dirichlet process, its conditional distribution is again a Dirichlet process. Moreover,

$$\mathbb{P}(\mathbf{S}_{L_m^{(n)}} = (s_1, \ldots, s_{K_{m(n)}}) | K_m^{(n)} = k, L_m^{(n)} = s, K_n = j) = \frac{\prod_{i=1}^k (s_i - 1)!}{|\mathfrak{s}(s, k)|}.$$

A further example of exchangeable Gibbs-type random partition for which closed-form expressions of sampling formulae are available is the generalized gamma distribution (6). The conditional EPPF of the corresponding random partition is given by

$$(27) \quad \begin{aligned} &\tilde{\Pi}_k^{(s)}(s_1, \ldots, s_k; m, n, j) \\ &= \frac{\sigma^k \sum_{i=0}^{n+m-1} \binom{n+m-1}{i} (-1)^i \beta^{i/\sigma} \Gamma(j + k - i/\sigma; \beta)}{\sum_{i=0}^s \mathscr{C}(s, i, \sigma) \sum_{l=0}^{n+m-1} \binom{n+m-1}{l} (-1)^l \beta^{l/\sigma} \Gamma(j + i - l/\sigma; \beta)} \\ &\quad \times \prod_{i=1}^k (1 - \sigma)_{s_i - 1} \end{aligned}$$

and all sampling distributions described in Section 2 can be derived in a straightforward way.

3.3. *Looking backward.* In this section we face the problem of determining the probability that certain specific observations, present in the basic sample, are not re-observed in the additional $m$-sample. This is tantamount to deriving the probability that the new observations belong either to new clusters or to specified "old" clusters.

Let $A_1, \ldots, A_j$ be the classes of $K_n = j$ sets into which the first $n$ observations, or $n$ integers $\{1, \ldots, n\}$, are clustered. Define $M_r^{(n,j)} := M_r^{(n,j)}(i_1, \ldots, i_r)$ to be, for any $(i_1, \ldots, i_r) \in \{1, \ldots, j\}^r$ such that $i_k \neq i_l$ for any $l \neq k$, the event which is true if and only if none of the $m$ observations belongs to any of the sets $A_i$ where $i \notin \{i_1, \ldots, i_r\}$. That is, $M_r^{(n,j)}$ is true if the $m$ new observations belong either to "new" clusters or to the specified "old" clusters $A_{i_1}, \ldots, A_{i_r}$. We are now interested in evaluating the probability of such an event. Obviously, one has $r \in \{1, \ldots, j\}$ and recall that of the $m$ new observations $m - s$ are the ones belonging to the "old" clusters. Correspondingly, we set $\mathbf{\Lambda}_r = (\Lambda_{i_1}, \ldots, \Lambda_{i_r})$ to be the vector of frequencies, that is, $\Lambda_{i_l} = \text{card}(\{n+1, \ldots, n+m\} \cap A_{i_l}) \geq 0$ for any $l = 1, \ldots, r$ and $\sum_{l=1}^{r} \Lambda_{i_l} = m - s$. Hence, it can be seen that

$$\mathbb{P}(K_n = j, \mathbf{N}_n = \mathbf{n}_j, L_m^{(n)} = s, K_m^{(n)} = k, \mathbf{S}_{L_n^{(m)}} = \mathbf{s}_{K_n^{(m)}}, \mathbf{\Lambda}_r = \mathbf{\lambda}_r) \quad (28)$$

$$= V_{n+m, j+k} \prod_{r=1}^{k} (1-\sigma)_{s_r - 1} \prod_{l=1}^{r} (1-\sigma)_{n_{i_l} + \lambda_{i_l} - 1} \prod_{l=r+1}^{j} (1-\sigma)_{n_{i_l} - 1}.$$

From (28) a number of interesting distributions can be derived. They typically provide information about the possibility of not re-observing certain "old" species in a subsequent "new" sample. The main result of the subsection we wish to state is the following:

PROPOSITION 4. *Given that the basic $n$-sample is partitioned into $K_n = j$ classes, $A_1, \ldots, A_j$, with frequencies $(n_1, \ldots, n_j)$, the probability that the observations from the subsequent $m$-sample contain either elements from $A_{i_1}, \ldots, A_{i_r}$, with $r \in \{1, \ldots, j\}$, or from new clusters is given by*

$$\mathbb{P}(M_r^{(n,j)} | K_n = j, \mathbf{N}_n = \mathbf{n}_j) \quad (29)$$

$$= \sum_{k=0}^{m} \frac{V_{n+m, j+k}}{V_{n,j}} \frac{\mathscr{C}(m, k; \sigma, r\sigma - \sum_{l=1}^{r} n_{i_l})}{\sigma^k}.$$

For the two-parameter Poisson–Dirichlet process, one has

$$\frac{V_{n+m, j+k}}{V_{n,j} \sigma^k} = \frac{(\theta+1)_{n-1} \prod_{i=1}^{j+k-1} (\theta + i\sigma)}{\sigma^k (\theta+1)_{n+m-1} \prod_{i=1}^{k-1} (\theta + i\sigma)}$$

$$= \frac{\prod_{i=0}^{k-1} (\theta + j\sigma + i\sigma)}{\sigma^k (\theta + n)_m} = \frac{((\theta + j\sigma)/\sigma)_k}{(\theta + n)_m}.$$

Hence, combining (29) with the definition of noncentral generalized factorial coefficient (38) in the Appendix, one has

$$
\mathbb{P}(M_r^{(n,j)}|K_n = j, \mathbf{N}_n = \mathbf{n}_j)
$$
$$
(30) \quad = \frac{1}{(\theta+n)_m} \sum_{k=0}^{m} \left(\frac{\theta+j\sigma}{\sigma}\right)_k \mathscr{C}\left(m, k; \sigma, r\sigma - \sum_{l=1}^{r} n_{i_l}\right)
$$
$$
= \frac{(\theta + (j-r)\sigma + \sum_{l=1}^{r} n_{i_l})_m}{(\theta+n)_m}.
$$

Such a simple expression provides the conditional probability that no integer in $\{n+1, \ldots, n+m\}$ will belong to any of the sets $A_i$, with $i \notin \{i_1, \ldots, i_r\}$, generated by $[n]$. In other terms, of the $j$ clusters associated to the (conditioning) partition of $[n]$, at most the $r$ clusters with indexes $i_1, \ldots, i_r$ do possibly contain integers from $\{n+1, \ldots, n+m\}$.

3.4. *The case $\sigma < 0$.* In the previous subsections we have focused on Gibbs random partitions with $\sigma \in [0, 1)$. See Definition 1. The nonnegativity of $\sigma$ ensures that (7) defines the probability distribution of an infinite exchangeable partition. On the other hand, when $\sigma < 0$ Lemma 8 in [9] entails that the weights $V = \{V_{n,k} : k = 1, \ldots, n; n \geq 1\}$ in (7) are mixtures of the weights $V(\nu) = \{V_{n,k}(\nu) : k = 1, \ldots, n; n \geq 1\}$, for $\nu = 1, 2, \ldots$, and $V_{n,k}(\nu) = |\sigma|^k \nu(\nu-1)\cdots(\nu-k+1)/(\nu|\sigma|)_n$. Thus, the $V_{n,k}(\nu)$'s correspond, for each $\nu = 1, 2, \ldots$, to the weights of Pitman's sampling formula (25) and imply that the number of distinct species in the population is $\nu$. See [9] for details. According to Theorem 12(i) in [9], $V$ arises as a mixture of the weights $V(1), \ldots, V(N^*)$, where $N^* \in \{1, 2, \ldots\} \cup \{\infty\}$ and one can, then, obtain the same results as stated in the present section, with the proviso that $K_n < N^* + 1$. If $N^* = \infty$, then no relevant change occurs. In particular, a slight modification of the proof allows one to recover the characterization of Corollary 2, with the conditional EPPF $\tilde{\Pi}_k^{(s)}(s_1, \ldots, s_k; m, n, j)$ being defined for any $j < N^* + 1$. Note that, in this case, the two-parameter Poisson–Dirichlet model coincides with symmetric Dirichlet distributions. See [29].

**4. Application to the analysis of EST data.** In this section we show how Gibbs priors can be applied in a straightforward way to the analysis of Expressed Sequence Tags (ESTs). ESTs are generated by partially sequencing randomly isolated gene transcripts that have been converted into cDNA. From their introduction in [1], ESTs have played an important role in the identification, discovery and characterization of organisms as they provide an attractive and efficient alternative to full genome sequencing. The resulting transcript sequences and their corresponding abundances are the main focus of interest providing the identification and level of expression of genes. Given a cDNA library and an initial sample of reads of size $n$, the main statistical issues to be faced are of predictive nature in the sense that various features of a possible additional sample of size $m$ are to be predicted. See, for example, [23, 24, 32]. Such features include, for instance: (i) the expected number of new genes meant as an estimate of the number of new unique genes to be detected in the additional EST survey; (ii) the expected number of genes which do not coincide with genes already present in the initial sample; (iii) the probability that certain specific genes, present in the basic sample, do not appear in the additional sample. Based on these estimates important decisions are to be taken. For instance, researchers have to decide: (i) whether to proceed with sequencing from a certain library; (ii) whether to carry out a "normalization" protocol (an expensive procedure which aims at making the frequencies of genes in the library more uniform); (iii) which libraries, among several ones concerning the same organism, are less redundant in the sense that they deliver more information from an additional sample.

The Bayesian nonparametric framework based on Gibbs-type random probability measures represents a natural, and at the same time powerful, approach for dealing with these kinds of problems since it conveys, in a statistically rigorous way, the information present in the initial sample into prediction. In particular, we focus on the two-parameter Poisson–Dirichlet process, which stands out for its mathematical tractability.

In order to illustrate the results of the previous section we first deal with a simple numerical example and then analyze some real EST data. The information provided by an EST data set sequenced from a cDNA library is summarized by the size of the sample $n$, the number of different cDNA fragments $j$, each of which represents a unique gene and their corresponding expression levels. Recalling the notation set in the Introduction, $M_{i,n}$ stands for the number of clusters of size $i$ with the initial $n$-sample: within the EST framework $M_{i,n}$ is now the number of genes with expression level $i$. For our purposes it is useful to convert the $M_{i,n}$'s into the $N_{i,n}$'s, the frequencies (or expression levels) of the various unique genes: hence, the sample information is given by $n$, $j$ and $(n_1, \ldots, n_j)$. We then assume the EST data are an exchangeable sequence with nonparametric prior given by the two-parameter Poisson–Dirichlet process. This implies that the clustering of the ESTs follows a two-parameter Poisson–Dirichlet random partition (5). Such a setup postulates the sequence of tags to be extendible to infinity: however, interest relies in computing estimates for $m$ up to the size of the library, which is always finite implying finiteness of all the estimates. In order to specify the prior parameters $\theta$ and $\sigma$ we resort to an empirical Bayes approach as in [20]. Hence, we fix $\sigma$ and $\theta$ so to maximize (5) corresponding to the observed sample $(j, n_1, \ldots, n_j)$, that is,

$$(31) \qquad (\hat{\sigma}, \hat{\theta}) = \arg\max_{(\sigma,\theta)} \frac{\prod_{i=1}^{j-1}(\theta + i\sigma)}{(\theta + 1)_{n-1}} \prod_{i=1}^{j} (1 - \sigma)_{n_i - 1}.$$

Given this, the model is completely specified and attention can be focused on predicting various features of a future sample of size $m$.

4.1. *Numerical example.* Here we compare predictions arising from two different basic samples both of size $n = 100$. The sample sequenced from library 1 is composed of $j = 59$ unique genes with $m_{1,100} = 40$, $m_{2,100} = 10$, $m_{3,100} = 4$, $m_{4,100} = 2$, $m_{5,100} = 2$, $m_{10,100} = 1$, whereas the sample sequenced from library 2 consists of $j = 37$ unique genes such that $m_{1,100} = 20$, $m_{2,100} = 5$, $m_{3,100} = 4$, $m_{4,100} = 3$, $m_{5,100} = 2$, $m_{6,100} = 1$, $m_{10,100} = 1$, $m_{20,100} = 1$. It is to be noted that the first one features a higher number of unique genes and the expression levels of the genes is remarkably lower. The average expression level, $n/j$, is 1.69 for the sample taken from library 1 and 2.7 for the sample sequenced from library 2. The parameters for Pitman's sampling formula are set according to (31), which yields $(\hat{\sigma}_1, \hat{\theta}_1) = (0.34, 33)$ and $(\hat{\sigma}_2, \hat{\theta}_2) = (0.26, 12)$ for the two cases. Furthermore, we consider an additional sample of size $m = 100$.

The expected number of new genes in the additional $m$-sample can be immediately derived from (14) and is given by

$$(32) \qquad \mathbb{E}[K_m^{(n)} | K_n = j] = \sum_{k=1}^{m} k \frac{(j + \theta/\sigma)_k}{(\theta + n)_m} \mathscr{C}(n, k; \sigma, -n + j\sigma).$$

In our case the estimator leads to predict 33 and 15 new unique genes, respectively. This is in accordance with the intuition, which leads to guess a higher number of new genes for library 1, since the basic sample featured 59 unique ones in contrast with 37 of library 2. A second quantity of interest is the expected number of genes, in the additional sample, which do not coincide with previously observed ones. Such an expression is given in (24) and it can be seen as a better measure of redundancy of the library since, in contrast to (32), it takes also the expression levels of the new genes into account. In our case $\mathbb{E}[L_m^{(n)} | K_n = j]$ yields 40 for library 1 and 19 for library 2. At first glance these estimates may seem low since one would expect the difference $\mathbb{E}[L_m^{(n)} | K_n = j] - \mathbb{E}[K_m^{(n)} | K_n = j]$ to be larger. However, it is reasonable that only a few new unique genes will have expression levels greater than 1: otherwise they would have been discovered already in the basic sample. Combining the two estimates one can obtain a plug-in estimator of the average expression level of the new unique genes in the additional sample as $A_m := \mathbb{E}[L_m^{(n)} | K_n = j] / \mathbb{E}[K_m^{(n)} | K_n = j]$, which in our case are equal to 1.21 for library 1 and 1.28 for library 2. If one is interested in the overall average expression level after $n + m = 200$ reads, the estimator

$$A_{n+m} := (n+m)/(j + \mathbb{E}[K_m^{(n)}|K_n = j]) \tag{33}$$

yields 2.17 for library 1 and 3.85 for library 2.

Another important aspect to look at is represented by the frequency configurations of the new unique genes in the additional sample. In particular, one is interested in establishing which type of configurations are more likely to appear. By the above considerations, it is clear that these will have a few numbers of unique genes with significant expression level (which have "escaped" being sequenced in the basic sample) and all the others with expression level 1. For detecting such a feature, we work conditionally on $K_m^{(n)}$ and on $L_m^{(n)}$ which leads to the following probability distribution for $\mathbf{S}_{L_m^{(n)}}$

$$\mathbb{P}(\mathbf{S}_s = (s_1, \ldots, s_k)|K_m^{(n)} = k, L_m^{(n)} = s, K_n = j) \tag{34}$$
$$= \frac{\sigma^k}{\mathscr{C}(s,k,\sigma)} \prod_{i=1}^{k}(1-\sigma)_{s_i-1}.$$

It is then reasonable to set $K_m^{(n)}$ equal to the expected number of new unique genes arising from (32) and $L_m^{(n)}$ equal to the expected number of genes which coincide with any of the newly observed genes, given in (24) with $\sigma = \hat{\sigma}$ and $\theta = \hat{\theta}$ as in (31). Denote these values by $\bar{k}_m$ and $\bar{s}_m$, respectively. Given this we consider the ratio of the distribution in (34) for two configurations $\mathbf{S}_{\bar{s}_m}^1$ and $\mathbf{S}_{\bar{k}_m}^2$ as an index for establishing which configuration is more likely to appear. From (34) one immediately obtains

$$I_m^{(n)}(\mathbf{S}_{\bar{s}_m}^1, \mathbf{S}_{\bar{k}_m}^2) := \frac{\prod_{i=1}^{\bar{s}_m}(1-\sigma)_{s_i^1-1}}{\prod_{i=1}^{\bar{k}_m}(1-\sigma)_{s_i^2-1}} \tag{35}$$

where obviously $\sum_{i=1}^{\bar{k}_m} s_i^r = \bar{s}_m$. Let us first consider library 1 and compare $\mathbf{S}_{40}^1$ given by 32 genes with expression level 1 and 1 gene with expression level 8 with $\mathbf{S}_{40}^2$ such that 26 genes are observed once and 7 twice. Then, $I_{100}^{(100)}(\mathbf{S}_{40}^1, \mathbf{S}_{40}^2) = 34346$, that is, the unbalanced configuration with only one gene having expression level 8 is 34346 times more likely than most balanced configuration. If we compare the first configuration with $\mathbf{S}_{40}^3$ given by 31 genes with expression level 1 and two genes with expression levels 4 and 5, respectively, then it appears that configuration 1 is "only" 60 times more likely. For library 2, things are quite different, even though the unbalanced configuration still predominates. By comparing $\mathbf{S}_{19}^1$ given by 14 genes with expression level 1 and 1 gene with expression level 5 with $\mathbf{S}_{19}^2$, where 11 genes are observed once and 4 genes twice, one has $I_{100}^{(100)}(\mathbf{S}_{19}^1, \mathbf{S}_{19}^2) = 44$. This means that the odds in favor of the unbalanced configuration with respect to the most balanced are "only" 44. By taking an intermediate configuration such as 13 genes observed once and two observed 2 and 4 times, respectively, the odds reduce to 5.

Finally, it is also worth looking backward, in the sense of determining probabilities that certain genes present in the initial sample will not be re-observed in the additional survey of size $m$. Since one is typically interested in probabilities concerning the most highly expressed genes, or the genes with expression level 1, it is useful to order the frequencies in the initial sample in increasing order and denote them by $n_{(1)}, n_{(2)}, \ldots, n_{(j)}$. Then, from (30), the probability of not re-observing the $j - r$ most highly expressed genes is given by

$$
(36) \qquad \frac{(\theta + (j-r)\sigma + \sum_{i=1}^{r} n_{(i)})_m}{(\theta + n)_m}.
$$

In order to avoid that probabilities take on too-low values, set $m = 10$. As for library 1, the probability of not observing the unique gene with expression level 10 is 0.482, whereas the probability of not observing the 40 genes with expression level 1 is 0.118. It is also worth noting that the probability of not observing certain specific 10 genes with expression level 1 (out of the 40 present in the initial sample) is given by 0.611. From this, one can see that it is more likely to re-observe a gene with expression level $a$ than $a$ genes with expression level 1: this appears to be a reasonable and, indeed, desirable feature for a model dealing with species prediction problems. As for library 2, one can, for instance, compute that the probability of not re-observing the unique gene with expression level 20 is 0.156, while the probability of not re-sequencing the 20 genes with expression level 1 is 0.257. Again, the probability attached to highly expressed genes is more than proportional with respect to genes with expression level 1. Finally, note that these probabilities are not directly comparable between libraries: this is due to the fact that library 1 exhibits a higher estimate of new genes to be discovered in the additional sample and also a higher number of observations which belong to these new clusters: consequently, it is natural that the probabilities of not re-observing certain genes are always higher for library 1.

4.2. *Genomic example.* Here we analyze a tomato-flower cDNA library from the Institute for Genomic Research Tomato Gene Index with library identifier T1526 [30]. This library was made from 0–3 mm buds of tomato flowers and was previously analyzed in [20, 23, 24] with reference to the determination of the discovery probability of further reads from the library.

TABLE 1
*Description of the* 10 *sub-samples of size* $n = 1000$ *and the true values of the quantities to be estimated on the remaining* $m = 1586$ *data: the second column indicates the number* $K_{1000} = j$ *of distinct genes in the sub-sample; columns* 3–6 *report the true values of* $(K_{1586}^{(1000)}|K_{1000})$, $(L_{1586}^{(1000)}|K_{1000})$, $A_{1586}$ *and* $A_{2586}$

| $N$ | $j$ | $K^{\text{true}}$ | $L^{\text{true}}$ | $A_{1586}^{\text{true}}$ | $A_{2586}^{\text{true}}$ |
|---|---|---|---|---|---|
| 1  | 825 | 1000 | 1166 | 1.166 | 1.417 |
| 2  | 816 | 1009 | 1142 | 1.132 | 1.417 |
| 3  | 806 | 1019 | 1151 | 1.130 | 1.417 |
| 4  | 834 | 991  | 1146 | 1.156 | 1.417 |
| 5  | 820 | 1005 | 1150 | 1.144 | 1.417 |
| 6  | 831 | 994  | 1145 | 1.152 | 1.417 |
| 7  | 819 | 1006 | 1150 | 1.149 | 1.417 |
| 8  | 813 | 1012 | 1130 | 1.117 | 1.417 |
| 9  | 812 | 1013 | 1135 | 1.120 | 1.417 |
| 10 | 830 | 995  | 1157 | 1.163 | 1.417 |

The initial sample consists of $n = 2586$ ESTs with $j = 1825$ unique genes. The tomato flower data set shows the following expression levels:

$$m_{i,2586} = 1434, 253, 71, 33, 11, 6, 2, 3, 1, 2, 2, 1, 1, 1, 2, 1, 1$$

with $i \in \{1, 2, \ldots, 14\} \cup \{16, 23, 27\}$, which means that we are observing 1434 genes which appear once, 253 genes which appear twice, etc. The average expression level of the basic sample is 1.417.

We first perform a cross-validation study for assessing the performance of the method. To this end 10 sub-samples of size 1000 have been drawn without replacement from the available 2586 EST sample. On the basis of each sub-sample, the corresponding values of $(\sigma, \theta)$ have been fixed according to (31). Then, we have computed the estimators for an additional sample of size $m = 1586$, which corresponds to the remaining observed data. In addition to the Bayes estimates $\mathbb{E}(K_m^{(n)}|K_n = j)$ and $\mathbb{E}(L_m^{(n)}|K_n = j)$, we also computed, using the distributions of $(K_m^{(n)}|K_n = j)$ recoverable from (10) and of $(L_m^{(n)}|K_n = j)$ given in (23), the 95% highest posterior density (HPD) intervals; these represent the Bayesian counterpart to frequentist confidence intervals. Finally, also the estimates for the average expression levels have been computed. Table 1 reports the true values corresponding to each sub-sample, whereas Table 2 displays the estimates with corresponding 95% HPD intervals.

By comparing Table 1 and 2, one sees that 9 times out of 10 the highest posterior density interval covers the true number of distinct genes present in the additional sample, whereas the true number of genes not coinciding with previously observed ones is always covered. The average prediction errors are 24.5 and 21.2 genes, respectively. The average error in the estimation of the expression level of the additional sample is 0.0026, whereas the average error of the estimates of the overall expression level is 0.019. Given the fact that prediction is carried out over an additional sample of size about 1.5 times the used sub-sample, such results appear completely satisfactory.

We now deal with the problem of predicting the outcomes of future sequencing and, as possible sizes of the additional sample, we consider $m \in \{250, 500, 750, 1000\}$. As for the prior specification of $(\sigma, \theta)$ the maximization in (31) leads to $(\hat{\sigma}, \hat{\theta}) = (0.612, 741)$. The corresponding estimates for the expected number of new genes (32) and for the number of genes which do not coincide with previously observed ones (24) are reported in Table 2 together with the corresponding 95% HPD intervals. The estimates of the average expression level of the new unique genes and of the average expression level for the whole sample of size $n + m$ are also reported.

It is worth noting that for this real data set the two estimates in the first two columns are extremely close, leading to an extremely low average expression level for the new unique genes. This can be explained by two facts: (i) the number of genes with expression level 1 is already very high in the basic sample ($m_{1,2586} = 1434$); (ii) the basic sample is large ($n = 2586$) and, hence, it is very unlikely that several highly expressed genes have not been sequenced. In such a case the frequency configuration of the additional sample is forced to be unbalanced and there is no need to compute (35) to state this. Just note that the most balanced configuration for $m = 250$ would be 136 genes with expression level 1 and 2 genes with level 2.

TABLE 2

*Predictions, based on the sub-samples, of the quantities of interest on the remaining $m = 1586$ data: columns 2–3 display the parameter specifications derived from (31); columns 4–7 report the Bayes estimates $\hat{K} := E[K_{1586}^{(1000)}|K_{1000} = j]$ and $\hat{L} := E[L_{1586}^{(1000)}|K_{1000} = j]$ with the corresponding 95% highest posterior density (HPD) intervals; columns 8–9 display the estimates for the average expression level in the additional sample and in the whole sample denoted by $\hat{A}_{1586}$ and $\hat{A}_{2586}$, respectively*

| $N$ | $\hat{\sigma}$ | $\hat{\theta}$ | $\hat{K}$ | HPD 95% | $\hat{L}$ | HPD 95% | $\hat{A}_{1586}$ | $\hat{A}_{2586}$ |
|---|---|---|---|---|---|---|---|---|
| 1  | 0.72 | 444 | 1017 | (969, 1067)  | 1140 | (1089, 1190) | 1.121 | 1.404 |
| 2  | 0.75 | 344 | 1012 | (963, 1065)  | 1128 | (1076, 1180) | 1.115 | 1.415 |
| 3  | 0.78 | 254 | 1009 | (959, 1063)  | 1116 | (1063, 1170) | 1.106 | 1.425 |
| 4  | 0.75 | 410 | 1049 | (1001, 1101) | 1165 | (1115, 1215) | 1.111 | 1.373 |
| 5  | 0.65 | 583 | 979  | (932, 1028)  | 1118 | (1068, 1168) | 1.142 | 1.437 |
| 6  | 0.73 | 446 | 1034 | (986, 1084)  | 1155 | (1104, 1204) | 1.117 | 1.387 |
| 7  | 0.72 | 420 | 1004 | (955, 1055)  | 1128 | (1077, 1179) | 1.124 | 1.419 |
| 8  | 0.72 | 397 | 992  | (943, 1043)  | 1115 | (1063, 1167) | 1.124 | 1.433 |
| 9  | 0.69 | 457 | 976  | (928, 1028)  | 1108 | (1056, 1159) | 1.135 | 1.446 |
| 10 | 0.72 | 466 | 1027 | (980, 1078)  | 1151 | (1100, 1200) | 1.121 | 1.393 |

TABLE 3
*Estimates arising from the two-parameter PD model with $(\hat{\sigma}, \hat{\theta}) = (0.612, 741)$ for sizes of the additional sample corresponding to $m \in \{250, 500, 750, 1000\}$. $\hat{K}$ and $\hat{L}$ denote the Bayes estimates $\mathbb{E}[K_m^{(2586)} | K_{2586} = 1825]$ and $\mathbb{E}[L_m^{(2586)} | K_{2586} = 1825]$, respectively, and are reported together with their 95% highest posterior density (HPD) intervals. Columns 6 and 7 display the estimated average expression levels*

| $m$ | $\hat{K}$ | 95% HPD | $\hat{L}$ | 95% HPD | $\hat{A}_m$ | $\hat{A}_{2586+m}$ |
|---|---|---|---|---|---|---|
| 250 | 138 | (122, 156) | 140 | (124, 155) | 1.014 | 1.445 |
| 500 | 272 | (249, 297) | 279 | (256, 302) | 1.026 | 1.471 |
| 750 | 402 | (373, 433) | 419 | (390, 448) | 1.042 | 1.498 |
| 1000 | 530 | (496, 566) | 558 | (523, 593) | 1.053 | 1.522 |

Another issue of interest is the determination of the probability of not re-observing certain particular genes in the additional $m$-sample. This can be achieved via the expressions in (30) and (36). With reference to the EST data set we are analyzing, the probability of not re-observing genes with expression level larger than 10, which correspond to 9 genes with frequencies $11, 11, 12, 13, 14, 16, 16, 23, 27$, is given by $0.656, 0.123, 0.016$ for $m = 10, 50, 100$, respectively. The probability of not observing the 71 genes with expression level 3 is given by $0.593, 0.075, 0.006$ for $m = 10, 50, 100$, respectively.

## APPENDIX

**A.1. Generalized factorial coefficients.** The results in the previous sections rely on the generalized factorial coefficients: here we provide a short account of their definitions and of formulae for their evaluation. For further details and pointers to the literature, the reader can refer to [5, 6]. See also [9]. For any $n \geq 1$ and $k = 0, \ldots, n$, the generalized factorial coefficient $\mathscr{C}(n, k; \sigma)$ coincides with the coefficient of the $k$th order factorial of $t$ in the expansion of the $n$th order generalized factorial of $t$ with scale parameter $\sigma \in \mathrm{R}$, that is,

$$(\sigma t)_n = \sum_{k=0}^{n} \mathscr{C}(n, k; \sigma)(t)_k.$$

In order to determine the distribution of the number of different species appearing in a sample of size $n$, that is, $K_n$, we have resorted to the following representation:

(37) $$\mathscr{C}(n, k; \sigma) = \frac{1}{k!} \sum_{j=0}^{k} (-1)^j \binom{k}{j} (-j\sigma)_n$$

with the proviso that $\mathscr{C}(0,0;\sigma) = 1$ and $\mathscr{C}(n,0;\sigma) = 0$ for all $n \geq 1$. It is to be noted that $\mathscr{C}$ slightly differs from the definition of generalized factorial coefficient $C(n,k;\sigma)$ as given, for example, in [5, 6]. Indeed, one has $\mathscr{C}(n,k;\sigma) = (-1)^{n-k} C(n,k;\sigma)$.

Besides $\mathscr{C}(n,k;\sigma)$ we consider another quantity $\mathscr{C}(n,k;\sigma,\gamma)$ which is known as noncentral generalized factorial coefficient. It is defined as the coefficient of the $k$th order factorial of $t$ in the expansion of the $n$th order noncentral generalized factorial of $t$, with scale parameter $\sigma$ and noncentrality parameter $\gamma$, that is,

$$(38) \qquad (\sigma t - \gamma)_n = \sum_{k=0}^{n} \mathscr{C}(n,k;\sigma,\gamma)(t)_k.$$

Note that in [5] the definition of noncentral generalized factorial coefficient designates a quantity $C(n,k;\sigma,\gamma) = (-1)^{n-k}\mathscr{C}(n,k;\sigma,\gamma)$. From (2.60) in [5] it is seen that it can be represented as

$$(39) \qquad \mathscr{C}(n,k;\sigma,\gamma) = \frac{1}{k!} \sum_{j=0}^{k} (-1)^j \binom{k}{j} (-\sigma j - \gamma)_n$$

and this can be usefully employed in order to evaluate the probability of discovering a new species. Moreover, from (2.56) in [5] it is possible to establish a connection between noncentral and central generalized factorial coefficients

$$(40) \qquad \mathscr{C}(n,k;\sigma,\gamma) = \sum_{s=k}^{n} \binom{n}{s} \mathscr{C}(s,k;\sigma)(-\gamma)_{n-s}.$$

Finally we briefly recall the relation to Stirling numbers. Indeed,

$$\lim_{\sigma \to 0} \frac{\mathscr{C}(n,k;\sigma)}{\sigma^k} = |s(n,k)|$$

where, as before, $|s(n,k)|$ is the signless Stirling number of the first kind. Moreover, one has

$$\lim_{\sigma \to 0} \frac{\mathscr{C}(n,k;\sigma,\gamma)}{\sigma^k} = \sum_{i=k}^{n} \binom{n}{i} |s(i,k)|(-\gamma)_{n-i}.$$

**A.2. Multivariate Chu–Vandermonde formula.** Here we present a multivariate version of the celebrated Chu–Vandermonde identity. In, for example, [5] the following version of the Chu–Vandermonde identity is presented:

$$(41) \qquad [a+b]_n = \sum_{r=0}^{n} \binom{n}{r} [a]_r [b]_{n-r}$$

for any $a$ and $b$ in R, where $[x]_n := x(x-1)\cdots(x-n+1)$ stands for the descending factorial. Since a multivariate version in terms of rising factorials seems not readily available in the literature we present it together with a proof.

LEMMA A.1. *For each $q, j \geq 1$, set $A_{j,q} = \{(q_1, \ldots, q_j) : q_i \geq 0, \sum_{i=1}^{j} q_i = q\}$. Then*

$$(42) \quad \sum_{(q_1,\ldots,q_j) \in A_{j,q}} \binom{q}{q_1 \cdots q_j} \prod_{i=1}^{j} (a_i)_{n_i + q_i - 1} = \left(n - j + \sum_{i=1}^{j} a_i\right)_q \prod_{i=1}^{j} (a_i)_{n_i - 1}$$

*where $(n_1, \ldots, n_j)$ is such that $n_i > 0$, for $i = 1, \ldots, j$ and $\sum_{i=1}^{j} n_i = n$.*

PROOF. Since $\frac{\Gamma(a+n)}{\Gamma(a)} = (a)_n = (-1)^n [-a]_n$, from identity (41) one deduces that

$$(43) \quad (a+b)_n = \sum_{r=0}^{n} \binom{n}{r} (a)_r (b)_{n-r}.$$

The proof now follows by inductive reasoning. Suppose the identity holds true for $j - 1$, that is,

$$\sum_{(q_1,\ldots,q_{j-1}) \in A_{j-1,q}} \frac{q!}{q_1! \cdots q_{j-2}! q_{j-1}!} (a_{j-1})_{n_{j-1} + q_{j-1} - 1} \prod_{i=1}^{j-2} (a_i)_{n_i + q_i - 1}$$

$$= \left(n - (j-1) + \sum_{i=1}^{j-1} a_i\right)_q \prod_{i=1}^{j-1} (a_i)_{n_i - 1},$$

and we show it holds for $j$ as well. Indeed, observe that

$$\sum_{(q_1,\ldots,q_j) \in A_{j,q}} \frac{q!}{q_1! \cdots q_{j-1}! q_j!} (a_j)_{n_j + q_j - 1} \prod_{i=1}^{j-1} (a_i)_{n_i + q_i - 1}$$

$$= \sum_{q_{j-1}=0}^{q} \frac{q!}{q_{j-1}!(q - q_{j-1})!} (a_{j-1})_{n_{j-1} + q_{j-1} - 1}$$

$$\times \sum_{(q_1,\ldots,q_{j-1}) \in A_{j-1,q-q_{j-1}}} \frac{(q - q_{j-1})!}{q_1! \cdots q_{j-2}! q_j!} (a_j)_{n_j + q_j - 1} \prod_{i=1}^{j-2} (a_i)_{n_i + q_i - 1}.$$

By the induction hypothesis, the second factor above equals

$$\left(n - n_{j-1} - (j-1) + \sum_{i=1}^{j} a_i - a_{j-1}\right)_{q - q_{j-1}} (a_j)_{n_j - 1} \prod_{i=1}^{j-2} (a_i)_{n_i - 1}.$$

Finally the proof is completed by virtue of (43), after noting that

$$(a_{j-1})_{n_{j-1}+q_{j-1}-1} = (a_{j-1})_{n_{j-1}-1}(n_{j-1}-1+a_{j-1})_{q_{j-1}}. \qquad \square$$

Lemma A.1 can also be proved by combining the last relation displayed in the proof with the definition of the multinomial-Dirichlet distribution, according to which $\sum_{(q_1,\ldots,q_j)\in A_{j,q}} \binom{q}{q_1\cdots q_j} \prod_{i=1}^{j}(a_i)_{n_i+q_i-1}/(n-j+\sum_{i=1}^{j}a_i)_q = 1$.

### A.3. Proofs.

PROOF OF PROPOSITION 1. An obvious point to start from is the following: if we have seen $n$ observations partitioned into $j$ distinct groups, then the conditional probability that the next $q \geq 1$ observations provide no new groups is $\prod_{l=1}^{q}(1-V_{n+l,j+1}/V_{n+l-1,j})$ which, using the recursive formula (8) for $V_{n,k}$, is given by

$$\mathbb{P}(K_q^{(n)}=0|K_n=j,\mathbf{N}_n=\mathbf{n}_j)$$
$$= \prod_{l=1}^{q}(n+l-1-j\sigma)\frac{V_{n+l,j}}{V_{n+l-1,j}} = (n-j\sigma)_q \frac{V_{n+q,j}}{V_{n,j}}.$$

On the other hand, suppose we have seen $n+q$ observations yielding $K_{n+q} = j$ groups. Then the conditional probability of obtaining $K_s^{(n+q)} = k$ new groups of sizes $s_1,\ldots,s_k$ from the next $s$ observations, where none of these coincides with the first $n+q$, is given by

$$\frac{V_{n+q+s,j+k}}{V_{n+q,j}} \prod_{i=1}^{k}(1-\sigma)_{s_i-1}$$

where $s_1+\cdots+s_k = s$. If we now set $q+s=m$, the conditional probability of obtaining new groups with respective frequencies $s_1,\ldots,s_k$ in the $m$ observations following on from $n$, given $K_n=j$, is, due to exchangeability, found by multiplying the two conditional probabilities above and including the $\binom{m}{s}$ term. Hence one achieves (9). Note that an alternative proof can be given by considering the joint distribution of $(K_n, \mathbf{N}_n, K_m^{(n)}, L_m^{(n)}, \mathbf{S}_{K_m^{(n)}}, \mathbf{\Lambda}_j)$ where $\Lambda_j = (\lambda_{1,m-L_m^{(n)}},\ldots,\lambda_{j,m-L_m^{(n)}})$ is the vector of nonnegative integers denoting the number of new observations in each of the $j$ groups into which the first $n$ observations are partitioned, and then by using Lemma A.1. $\square$

PROOF OF PROPOSITION 2. The proof works by induction. Let us first note that for any $m \geq 1$ one has

$$L_{m+1}^{(n)} = L_m^{(n)} + H_{n,m}$$

where $H_{n,m} = \mathbb{I}_{\mathcal{X}_n^c}(X_{n+m+1})$ and $\mathcal{X}_n = \{X_1, \ldots, X_n\}$. Let us first fix $m = 1$ and determine

$$E[L_2^{(n)}|K_n = j] = E[L_1^{(n)}|K_n = j] + E[H_{n,1}|K_n = j].$$

The first summand is clearly equal to $V_{n+1,j+1}/V_{n,j}$. As for the second summand, one can use the assumption of exchangeability which yields

$$E[H_{n,1}|K_n = j] = E[\mathbb{I}_{\mathcal{X}_n^c}(X_{n+2})|K_n = j]$$
$$= E[\mathbb{I}_{\mathcal{X}_n^c}(X_{n+1})|K_n = j]$$
$$= \frac{V_{n+1,j+1}}{V_{n,j}}.$$

Hence, (16) holds true for $m = 2$. Now, suppose (16) is valid for $m$ and let us show this implies it is still true for $m + 1$. This means we shall determine

$$E[L_{m+1}^{(n)}|K_n = j] = E[L_m^{(n)}|K_n = j] + E[H_{n,m}|K_n = j].$$

By assumption, $E[L_m^{(n)}|K_n = j] = mV_{n+1,j+1}/V_{n,j}$. Moreover, exchangeability again entails that the second summand above is $V_{n+1,j+1}/V_{n,j}$ and the conclusion follows. $\square$

PROOF OF PROPOSITION 3. This is straightforward and follows from taking the ratio between (9) and (11) in Corollary 1. $\square$

PROOF OF COROLLARY 2. Let

$$f_{j,\sigma}(s,k) := \frac{V_{n+m,j+k}}{\sum_{i=0}^s V_{n+m,j+i}\sigma^{-i}\mathscr{C}(s,i;\sigma)}$$

which, by assumption, does not depend on $n$ and $m$. Then, if $s = 2$ and $k = 2$

$$\frac{V_{n+m,(j-2)+2}}{\sum_{i=0}^2 V_{n+m,(j-2)+i}\sigma^{-i}\mathscr{C}(2,i;\sigma)} = f_{j-2,\sigma}(2,2)$$

and, with $s = 2$ $k = 1$

$$\frac{V_{n+m,(j-2)+1}}{\sum_{i=0}^2 V_{n+m,(j-2)+i}\sigma^{-i}\mathscr{C}(2,i;\sigma)} = f_{j-2,\sigma}(2,1).$$

If we, now, consider the ratio of these two expressions, we obtain the identity

$$\frac{V_{n+m,j}}{V_{n+m,j-1}} = \frac{f_{j-2,\sigma}(2,2)}{f_{j-2,\sigma}(2,1)}.$$

From this, one sees that

$$\frac{V_{n+m,j}}{V_{n+m,2}} = \prod_{i=3}^j \frac{V_{n+m,i}}{V_{n+m,i-1}} = \prod_{i=1}^{j-2} \frac{f_{i,\sigma}(2,2)}{f_{i,\sigma}(2,1)}$$

so that $V_{n+m,j} = V_{n+m,2} \prod_{i=1}^{j-2}(f_{i,\sigma}(2,2)/f_{i,\sigma}(2,1)) = g_1(n+m)g_2(j)$ for some functions $g_1$ and $g_2$. By a result of Kerov ([18], Theorem 7.1) (see also [9]) this entails that the weights $V_{n,k}$ are those from the two-parameter Poisson–Dirichlet process. □

PROOF OF PROPOSITION 4. Consider the expression displayed in (28) and sum with respect to all vector $\mathbf{s}_k$ in $\Delta_{s,k}$ to obtain

$$
\begin{aligned}
(44) \quad & \mathbb{P}(K_n = j, \mathbf{N}_n = \mathbf{n}_j, L_m^{(n)} = s, K_m^{(n)} = k, \mathbf{\Lambda}_r = \boldsymbol{\lambda}_r) \\
& = V_{n+m,j+k} \frac{\mathscr{C}(s,k,\sigma)}{\sigma^k} \prod_{l=1}^{r}(1-\sigma)_{n_{i_l}+\lambda_{i_l}-1} \prod_{l=r+1}^{j}(1-\sigma)_{n_{i_l}-1}.
\end{aligned}
$$

Now exploit Lemma 1 in order to integrate out $\mathbf{\Lambda}_r$ thus obtaining

$$
\begin{aligned}
(45) \quad & \mathbb{P}(K_n = j, \mathbf{N}_n = \mathbf{n}_j, L_m^{(n)} = s, K_m^{(n)} = k, M_r^{(n)}) \\
& = V_{n+m,j+k} \frac{\mathscr{C}(s,k,\sigma)}{\sigma^k} \left(\sum_{l=1}^{r} n_{i_l} - r\sigma\right)_{m-s} \prod_{i=1}^{j}(1-\sigma)_{n_i-1}.
\end{aligned}
$$

Finally, integrating out $L_m^{(n)}$ and $K_m^{(n)}$ and summing over $k = 0, \ldots, m$ one has

$$
\begin{aligned}
(46) \quad & \mathbb{P}(K_n = j, \mathbf{N}_n = \mathbf{n}_j, M_r^{(n)}) \\
& = \prod_{i=1}^{j}(1-\sigma)_{n_i-1} \sum_{k=0}^{m} \frac{V_{n+m,j+k}}{\sigma^k} \mathscr{C}\left(m,k;\sigma,r\sigma - \sum_{l=1}^{r} n_{i_l}\right).
\end{aligned}
$$

Hence, the ratio of (46) over the EPPF $\Pi_j^{(n)}(n_1, \ldots, n_j)$ yields the result in (29). □

**Acknowledgments.** The authors are grateful to two anonymous referees for their valuable comments and suggestions. Special thanks also to Ramsés H. Mena for helpful advice.

A. LIJOI  
DIPARTIMENTO DI ECONOMIA POLITICA  
 E METODI QUANTITATIVI  
UNIVERSITÀ DEGLI STUDI DI PAVIA  
VIA SAN FELICE 5  
27100 PAVIA  
ITALY  
E-MAIL: lijoi@unipv.it

I. PRÜNSTER  
DIPARTIMENTO DI STATISTICA  
 E MATEMATICA APPLICATA  
UNIVERSITÀ DEGLI STUDI DI TORINO  
PIAZZA ARBARELLO 8  
10122 TORINO  
ITALY  
E-MAIL: igor@econ.unito.it

S. G. WALKER  
INSTITUTE OF MATHEMATICS, STATISTICS  
 AND ACTUARIAL SCIENCE  
UNIVERSITY OF KENT  
KENT CT2 7NZ  
UNITED KINGDOM  
E-MAIL: S.G.Walker@kent.ac.uk